\documentclass[11pt]{article}
\usepackage{amssymb,amsmath}
\usepackage{graphicx}
\usepackage[latin1]{inputenc}
\newtheorem{dfn}{Definition}[section]

\newtheorem{thm}{Theorem}[section]
\newtheorem{thm*}{Theorem}
\newtheorem{cor}{Corollary}[section]
\newtheorem{prop}{Proposition}[section]

\newtheorem{ex}{Example}[section]
\newcommand{\Pf}{{\em Proof}. }
\newcommand{\EPf}{\begin{flushright} $\Box$ \end{flushright}}

\newcommand{\ca}{{\rm ca}}

\begin{document}
\title{Rotation surfaces of constant Gaussian curvature  as Riemannian approximation scheme in sub-Riemannian Heisenberg space $\mathbb H^1$}
\author{ 
       Jos\'e M. M. Veloso
         \\ Instituto de Ci\^encias Exatas e Naturais
         \\ Universidade Federal do Par\'a
                  \\ veloso@ufpa.br
                  \\https://orcid.org/0000-0002-8969-7320 }

\maketitle
\begin{abstract}
We verify if  Gausssian curvature of surfaces and normal curvature of curves in surfaces introduced by Diniz-Veloso \cite{DV} and by Balogh-Tyson-Vecchi \cite{BTV} to prove Gauss-Bonnet theorems in Heisenberg space $\mathbb H^1$ are equal.    The authors in \cite{BTV} utilize a limit of Gaussian and normal curvatures defined in the Riemannian approximations scheme $(\mathbb R^3,g_L)$ introduced by Gromov to study sub-Riemannian spaces. They show that these limits exist (unlike the limit of Riemannian surface  area form or length form), and they obtain Gauss-Bonnet theorem in $\mathbb H^1$ as limit of Gauss-Bonnet theorems in $(\mathbb R^3,g_L)$ when $L$ goes to infinity. The approach in \cite{DV} uses an analogous of Gauss map defined in the unitary cylinder and they show that the curvature defined by the Gauss map  corresponds to the curvature of covariant derivative which is null on left invariant vector fields of $\mathbb H^1$. The proof of Gauss-Bonnet theorem follows as the classical one. 
Curvatures  do not coincide. For comparison sake we apply the same formalism of \cite{DV} to get the curvatures of \cite{BTV}. With the obtained formulas, it is possible to prove the Gauss-Bonnet theorem in \cite{BTV} as a straightforward application of Stokes theorem. 

To exemplify the Gaussian curvature of \cite{BTV},  we calculate  the  rotation surfaces of  constant curvature in  $\mathbb H^1$, which are  only of three types. 


\end{abstract}


\section{Introduction} In \cite{DV}  Gaussian curvature for non-horizontal surfaces in sub-Riemannian Heisenberg space $\mathbb H^1$ was defined and  a Gauss-Bonnet Theorem was proved. The definition was analogous to Gauss curvature of surfaces in $\mathbb R^3$ with particular normal to surface and Hausdorff measure of area. The image of Gauss map was in the cylinder of radius one.

In a recent paper \cite{BTV}, a version of Gauss-Bonnet theorem is given through a limit of Riemannian approximations. The objective of this paper is to compare curvatures defined in \cite{BTV} and \cite{DV}, and show  that these are not coincident. Furthermore, we show that the limit $K^\infty$ and $k_n$ of Riemannian Gaussian curvatures $K^L$ and normal curvatures $k_n^L$ depend only on one of the two functions that define the geometry of the surface. Also we obtain that $K^\infty d\sigma=d(k_nds)$ and  it is possible to get Gauss-Bonnet theorem of \cite{BTV} applying Stokes theorem without taking limit.

The space $\mathbb H^1$ is a Lie group. We define a distribution $D$ generated by the left invariant vector fields $e_1=\frac{\partial}{\partial x}-\frac 1 2y\frac{\partial}{\partial z}$ and $e_2=\frac{\partial}{\partial y}+\frac 1 2x\frac{\partial}{\partial z}$ on $\mathbb H^1$. Introduce a scalar product in $D$ such that $e_1,e_2$ are orthonormal. We complete $e_1,e_2$ to a basis of left invariant vector fields in $\mathbb H^1$ by introducing $e_3=[e_1,e_2]=\frac{\partial}{\partial z}$. 
In \cite{BTV}, \cite{CDPT}  they consider the family of $g^L$ metrics such that $e_1$, $e_2$, $e_3^L=e_3/\sqrt{L}$ is an orthonormal basis,
and the Levi-Civitta connection $\overline\nabla^L$ on $(\mathbb H^1,g^L)$.

At points of a surface $S$, where the distribution $D$ does not coincide with $TS$, the intersection $D\cap TS$ has dimension 1, and we obtain a direction that we call \emph{characteristic} at this point of $S$. We suppose that all points of  surface $S$ have this property. 
The vector field  \emph{normal horizontal} $f_1$ is a unitary vector field in $D$ orthogonal to the characteristic direction which we suppose to be globally defined. We define $f^1$ by $f^1(f_1)=1$ and $f^1(TS)=0$. We denote by $f_2$ a unitary vector field in $D\cap TS$ and complete a basis of $TS$ taking $f_3=e_3+Af_1$. If $\alpha$ is the angle between $e_1$ and $f_1$, then $f_1=\cos\alpha\, e_1+\sin\alpha\, e_2$ and $f_2=-\sin\alpha\, e_1+\cos\alpha\, e_2$.
An orthonormal basis of $TS$ in $(\mathbb R^3,g^L)$ is given by $X_2^L=f_2$, $X_3^L=\frac{1}{\sqrt{L+A^2}}f_3$. The normal vector in $g^L$ to $S$ is 
$$X_1^L=\frac{\sqrt{L}}{\sqrt{L+A^2}}f_1-\frac{A}{\sqrt{L+A^2}}e_3^L.$$ 
The curvature of $S$ in the metric $g^L$ is
$$K^L=\frac{L}{(L+A^2)^2} d\alpha\wedge dA(f_3,f_2)-\frac{L^2}{(L+A^2)^2}dA(f_2)-\frac{L}{L+A^2}A^2,$$
and 
$$K^\infty=\lim_{L\rightarrow \infty}K_L=-dA(f_2)-A^2.$$
This formula shows $K^\infty$ depends only on $A$.

We briefly introduce the curvature $K$ of $S$ as  in \cite{DV}.  We consider  the dual horizontal normal $f^1=\cos\alpha e^1+\sin\alpha e^2-Ae^3$ as a Gauss map 
$$
g:p\in S\rightarrow (\cos\alpha,\sin\alpha, -A)\in C,
$$
where $C=\{(x,y,z):x^2+y^2=1\},$
and
$$
K(p)=\lim_{U\rightarrow\{p\}}\frac{\int_{g(U)}i(\tilde f_1)dV}{\int_{U}i(f_1)dV}
$$
is the Gaussian curvature of surface $S$ at point $p$, with $dV=f^1\wedge f^2\wedge f^3$ and $\tilde f_1$ the horizontal normal to $C$.
A simple calculation shows that 
$$K=d\alpha\wedge dA(f_3,f_2).$$
Observe that the curvature $K$ appears in the expression of $K^L$ as a term that goes to $0$ when $L\rightarrow \infty$.

We finish this work by calculating the surfaces invariant by rotations around $z$-axis  with $K^\infty$ constant.
A surface invariant by rotations is foliated by the  horizontal curves  which are tangents to $f_2$  (they are unique up to rotations),.  We write the ordinary differential equations which give these surfaces in terms of horizontal curves. If the rotation surface is obtained by rotating the curve $(r(t),0,c(t))$ around the $z$-axis, then by a parameterization satisfying
\begin{equation}\label{edoteta1}
\theta'=\frac{\sqrt{1- (r')^2}}{r}
\end{equation}
we obtain the other equations 
\begin{equation}\label{edor}
K^\infty+\frac{d^2}{dt^2}(\ln r^2)+\frac{d}{dt}(\ln r^2)=0,
\end{equation}
\begin{equation}\label{edoc}
c'=\frac{1}{2}r \sqrt{1-(r')^2}.
\end{equation}

In case of constant Gaussian curvature $K^\infty$, we integrate these equations and give the graphics of these surfaces  which are solutions of (\ref{edoteta1}), (\ref{edor}), (\ref{edoc}).  There exists only three classes of rotation surfaces with $K^\infty$ constant, which is different from  the numerous  rotation surfaces of constant curvature in Euclidean space $\mathbb R^3$ or those in $\mathbb H^1$ with $K$ constant. 



For all the surface graphics we show  the horizontal curves inside them.

For more details on these topics see  
\cite{CDPT}, \cite{DSV}.

This work was supported by CAPES.

\section{The Heisenberg group}
We denote by $\mathbb H^1$ the Heisenberg nilpotent Lie group whose manifold is $\mathbb R^3$, with Lie algebra ${H}^1=V_1\oplus V_2$, $\dim V_1=2$, $\dim V_2=1$ and  $[V_1,V_1]=V_2$, $[V_1,V_2]=0$. Since   $\mathbb H^1$ is nilpotent, the exponential map $\exp:H^1\rightarrow \mathbb H^1$ is a diffeomorphism. Let  $e_1,e_2$  be a basis of $V_1$ and $e_3=[e_1,e_2]\in V_2$. Applying the Baker-Campbell-Hausdorff formula we have
$$\exp^{-1}(\exp(X)\exp(Y))=X+Y+\frac{1}{2}[X,Y],$$
where $X=x_1e_1+y_1e_2+z_1e_3$ and $Y=x_2e_1+y_2e_2+z_2e_3$.  As $[e_1,e_2]=e_3$ we get 
$$ X+Y+\frac{1}{2}[X,Y]=(x_1+x_2)e_1+(y_1+y_2)e_2+(z_1+z_2+\frac{1}{2}(x_1y_2-x_2y_1))e_3. $$ 
We identify $\mathbb H^1$ with $\mathbb R^3$ by identifying $(x,y,z)$ with $\exp(xe_1+ye_2+ze_3)$, and this identification is known as  canonical coordinates of first kind or exponential coordinates. In these coordinates the group operation is
$$(x_1,y_1,z_1)(x_2,y_2,z_2)=(x_1+x_2,y_1+y_2,z_1+z_2+\frac{1}{2}(x_1y_2-x_2y_1)),$$
the exponential is 
$$\exp(xe_1+ye_2+ze_3)=(x,y,z),$$
and the left invariant vector fields $e_1,e_2,e_3$ are given by
$$\begin{array}{lcl}e_1&=&\frac{\partial}{\partial x}-\frac 1 2y\frac{\partial}{\partial z},\vspace{.2cm}\\
e_2&=&\frac{\partial}{\partial y}+\frac 1 2x\frac{\partial}{\partial z},\vspace{.2cm}\\
$$e_3&=&\frac{\partial}{\partial z},
\end{array}
$$
with brackets 
$[e_1,e_2]=e_3,$ and $[e_3,e_1]=[e_3,e_2]=0$.
 The dual basis is
$$\begin{array}{lcl}
e^1&=&dx,\vspace{.2cm}\\
e^2&=&dy,\vspace{.2cm}\\
e^3&=&dz+\frac 1 2(ydx-xdy),
\end{array}
$$
with $de^3=-e^1\wedge e^2$, $de^1=de^2=0$.

Let  $D\subset T\mathbb H^1$ be the two dimensional distribution generated by the vector fields $e_1,e_2$, so that $D$ is the null space of $e^3$. On $D$ we define a scalar product $<\, ,\, >$ such that $\{e_1,e_2\}$ is an orthonormal basis of $D$. 

The \emph{element of volume} $dV$ in $\mathbb H^1$ is $dV=e^1\wedge e^2\wedge e^3=dx\wedge dy\wedge dz.$ 

\section{The approximation by scalar product $g_L$}

Consider the $g_L$ metrics where $e_1$, $e_2$, $e_3^L=e_3/\sqrt{L}$ is an orthonormal basis. The dual basis is $e^1$, $e^2$, $e^3_L=\sqrt{L}e^3$. Then the Carnot-Caratheodory metric space $\mathbb H^1$ is the limit in the sense of Gromov-Hausdorff of Riemannian metric spaces $(R^3,d_L)$, when $L\rightarrow \infty$.  We consider the Levi-Civitta connection $\overline\nabla^L$ in $(\mathbb H^1,g_L)$ . A straightforward calculation shows that (see \cite{CDPT}) :
\begin{equation}\label{overnabla}
\begin{array}{rcl}
\overline\nabla^Le_1&=&\frac{\sqrt{L}}{2}(-e^2\otimes e_3^L- e_L^3\otimes e_2)\\
\overline\nabla^Le_2&=&\frac{\sqrt{L}}{2}(e^1\otimes e_3^L+ e_L^3\otimes e_1)\\
\overline\nabla^Le_3^L&=&\frac{\sqrt{L}}{2}(-e^1\otimes e_2+ e^2\otimes e_1).
\end{array}
\end{equation}
The curvatute tensor
$\overline R_{L,ijkl}=\overline R_L(e_i,e_j,e_k,e_l)$ is given by
$$\overline R_{L,1212}=\overline R_{L,2121}=-\overline R_{L,1221}=-\overline R_{L,2112}=\frac{3L}{4},$$
$$\overline R_{L,1313}=\overline R_{L,3131}=\overline R_{L,2323}=\overline R_{L,3232}=-\overline R_{L,1331}=-\overline R_{L,3113}=-\overline R_{L,2332}=-\overline R_{L,3223}=\frac{-L}{4},$$
and $\overline R_{L,ijkl}=0$ otherwise.

\section{Surfaces in $\mathbb H^1$}
Suppose that $S$ is an oriented differentiable two dimensional manifold in $\mathbb H^1$. We get that $\dim(D\cap TS)\geq 1$, and as $de^3=e^1\wedge e^2$, the set where $\dim(D\cap TS)=2$ has empty interior. We denote by $\Sigma=\{x\in S:  \dim(D_x\cap T_xS)=2\}$ and by $S'=S-\Sigma$. The set $S'$ is open in $S$. In what follows we will suppose $\Sigma=\emptyset$, so  $S=S'$. With this hypothesis the one dimensional vector subbundle $D\cap TS$  is well defined on $S$. Suppose $U\subset S$ is an open set such that we can define a unitary vector field $f_2$ with values in $D\cap TS$, so that $<f_2,f_2>=1$. Suppose $f_2=xe_1+ye_2$.
\begin{dfn}The unitary vector field $f_1\in \underline{D}$ defined by
$$f_1=ye_1-xe_2$$
is the \emph{horizontal normal} to $S$.
\end{dfn}
Then we can define $f^1\in \underline{(T\mathbb H^1)^*|_S}$ by $f^1(f_1)=1$ and $f^1(TS)=0$. We call $f^1$ the \emph{horizontal conormal} to $S$. 
If 
$$f_3=e_3-f^1(e_3)f_1,$$
then $\{f_2,f_3\}$ is a \emph{special} basis of $TS$ on the open set $U$. If we write
$$f_1=\cos\alpha\, e_1+\sin\alpha\, e_2,$$
for some real function $\alpha$ on $U$, reducing $U$ if necessary, then
$$f_2=-\sin\alpha\, e_1+\cos\alpha\, e_2,$$
and if we denote by $A=-f^1(e_3)$,
$$f_3=e_3+Af_1.$$
The dual basis of $(T\mathbb H^1)^*$ on $S$ is 
$$\begin{array}{lcl}
f^1&=&\cos\alpha\, e^1+\sin\alpha\, e^2-Ae^3,\\
f^2&=&-\sin\alpha\, e^1+\cos\alpha\, e^2,\\
f^3&=&e^3.
\end{array}$$
The inverse relations are
$$
\begin{array}{lcl}
e^3&=&f^3,\\
e^1&=&\cos\alpha\,f^1-\sin\alpha\, f^2+A\cos\alpha\, f^3,\\
e^2&=&\sin\alpha\,f^1+\cos\alpha\, f^2+A\sin\alpha\, f^3.
\end{array}
$$
As $f^1=0$ on $S$, we get
$$
\begin{array}{rcl}
df^2&=&-A\,d\alpha\wedge f^3,\\
df^3&=&Af^2\wedge f^3,\\
0&=&(d\alpha+A^2f^3)\wedge f^2-dA\wedge f^3.
\end{array}
$$
From this last relation we obtain $$d\alpha(f_3)=-(dA(f_2)+A^2).$$
We have on $TS$
\begin{equation}\label{dAf3}
\begin{array}{rcl}
d(Af^3)&=&dA\wedge f^3+Ade^3\\
&=&i(f_2)dAf^2\wedge f^3-Ae^1\wedge e^2\\
&=&i(f_2)dAf^2\wedge f^3-A(-\sin\alpha f^2+A\cos\alpha f^3)\wedge(\cos\alpha f^2+A\sin\alpha f^3)\\
&=&(i(f_2)dA+A ^2)f^2\wedge f^3
\end{array}
\end{equation}

\section{The orthonormal basis of $g_L$}
An orthonormal basis of $TS$ in $(\mathbb R^3,g_L)$ is given by \begin{equation}\label{X2X3}
X_2^L=f_2 \mbox{ and }X_3^L=\frac{1}{\sqrt{L+A^2}}f_3.
\end{equation} 
The normal vector in $g_L$ to $S$ is 
\begin{equation}\label{X1}
X_1^L=\frac{\sqrt{L}}{\sqrt{L+A^2}}f_1-\frac{A}{\sqrt{L+A^2}}e_3^L
=\frac{\sqrt{L+A^2}}{\sqrt{L}}f_1-\frac{A}{\sqrt{L+A^2}\sqrt{L}}f_3.
\end{equation} 
If we write $\cos\beta=\frac{\sqrt{L}}{\sqrt{L+A^2}}$ and $\sin\beta=\frac{A}{\sqrt{L+A^2}}$, then our orthonormal basis is 
$$
\begin{array}{rcl}
X_1^L&=&\cos\beta \cos\alpha e_1+\cos\beta\sin\alpha e_2-\sin\beta e_3^L\\
X_2^L&=&-\sin\alpha e_1+\cos\alpha e_2\\
X_3^L&=&\sin\beta\cos\alpha e_1+\sin\beta\sin\alpha e_2+\cos\beta e_3^L
\end{array}
$$
As $d\sin\beta=\cos\beta d\beta$, we get
$$
d\beta=\frac{\sqrt{L}}{L+A^2}dA.
$$

\section{The projection $\nabla^L$ of $\overline\nabla^L$ on $TS$}
The connection $\nabla^L$ on $TS$ is defined by
$$\nabla^L_XY=\overline\nabla^L_XY-<\overline\nabla^L_XY,X_1^L>X_1^L,$$
for $X,Y$ sections of $TS$.
We have 
$$\nabla^L X^L_2=<\overline\nabla^LX^L_2,X^L_3>X_3^L.$$
Taking into account (\ref{overnabla}) we get
$$\overline\nabla^LX^L_2=(-d\alpha+\frac{\sqrt{L}}{2}e^3_L)\otimes(\cos\alpha e_1+\sin\alpha e_2)+\frac{\sqrt{L}}{2}(\cos\alpha e^1+\sin\alpha e^2)\otimes e_3^L$$
and
$$
\begin{array}{rcl}
e^1&=&\cos\beta \cos\alpha X^1_L-\sin\alpha X^2_L+\sin\beta\cos\alpha X^3_L\\
e^2&=&\cos\beta\sin\alpha X^1_L+\cos\alpha X_L^2+\sin\beta\sin\alpha X^3_L\\
e^3_L&=&-\sin\beta X^1_L+\cos\beta X_L^3
\end{array}
$$
we get
$$
\begin{array}{rcl}
<\overline\nabla^LX^L_2,X^L_3>&=&(-d\alpha+\frac{\sqrt{L}}{2}e^3_L)\sin\beta+\frac{\sqrt{L}}{2}(\cos\alpha e^1+\sin\alpha e^2)\cos\beta\\
&=&-d\alpha\sin\beta+\frac{\sqrt{L}}{2}\cos(2\beta)X^1_L+\frac{\sqrt{L}}{2}\sin(2\beta)X^3_L
\end{array}
$$
As $X^1_L$ is null on $TS$ we get
$$
\nabla X_2^L=(-d\alpha\sin\beta+\frac{\sqrt{L}}{2}\sin(2\beta)X^3_L)\otimes X_3^L.
$$
In the same way as $<\overline\nabla^LX^L_3,X^L_2>=-<\overline\nabla^LX^L_2,X^L_3>$, we obtain
$$\nabla X_3^L=(d\alpha\sin\beta-\frac{\sqrt{L}}{2}\sin(2\beta)X^3_L)\otimes X_2^L.$$
Now 
$$
\begin{array}{rcl}
X^3_L([X^L_2,X^L_3])&=&\frac{1}{\cos\beta}e^3_L([X^L_2,X^L_3])=\\
&=&\frac{1}{\cos\beta}e^3_L\left((-\sin\beta\sqrt{L}+\sin\beta\sin\alpha d\beta(e_1)-\sin\beta\cos\alpha d\beta(e_2))e_3^L\right)\\
&=&\frac{\sin\beta}{\cos\beta}(-\sqrt{L}-d\beta(X_2^L))\
\end{array}
$$

\section{The limit of curvatures $K^L$ of $S$}\label{curvgauss}

Now we will calculate the Gaussian curvature $K^L=<R^L(X_2^L,X_3^L)X^L_3,X_2^L>$, where $R^L(X,Y)Z=\nabla_X\nabla_YZ-\nabla_Y\nabla_XZ-\nabla_{[X,Y]}Z$. Therefore
$$
<\nabla_{X_2^L}\nabla_{X_3^L}X_3^L,X_2^L>=X_2^L(d\alpha(X_3^L)\sin\beta-\frac{\sqrt{L}}{2}\sin(2\beta))
$$
$$
<\nabla_{X^L_3}\nabla_{X_2^L}X_3^L,X_2^L>=X_3^L(d\alpha(X_2^L)\sin\beta)
$$
$$
\begin{array}{rcl}
<\nabla_{[X^L_2,X^L_3]}X_3^L,X_2^L>&=&d\alpha([X^L_2,X^L_3])\sin\beta-\frac{\sqrt{L}}{2}\sin(2\beta)X^3_L([X^L_2,X^L_3])\\
&=&d\alpha([X^L_2,X^L_3])\sin\beta-\frac{\sqrt{L}}{2}\sin(2\beta)(\frac{\sin\beta}{\cos\beta}(-\sqrt{L}- d\beta(X_2^L)))\
\end{array}
$$
Then
$$
\begin{array}{rcl}
K^L&=&X_2^L(d\alpha(X_3^L)\sin\beta-\frac{\sqrt{L}}{2}\sin(2\beta))-X_3^L(d\alpha(X_2^L)\sin\beta)-d\alpha([X^L_2,X^L_3])\sin\beta\\
&&+\frac{\sqrt{L}}{2}\sin(2\beta)(\frac{\sin\beta}{\cos\beta}(-\sqrt{L}- d\beta(X_2^L)))\\
&=&d\alpha(X_3^L)\cos\beta d\beta(X^L_2)-\frac{\sqrt{L}}{2}2\cos(2 \beta)d\beta(X^L_2)-d\alpha(X_2^L)\cos\beta d\beta(X_3^L)\\
&&+\sqrt{L}\sin^2(\beta)(-\sqrt{L}- d\beta(X_2^L))\\
&=&\cos\beta d\alpha\wedge d\beta(X_3^L,X_2^L)-\sqrt{L}\cos^2\beta d\beta(X_2^L)-L\sin^2\beta\\
&=&\frac{\sqrt{L}}{\sqrt{L+A^2}} d\alpha\wedge \frac{\sqrt{L}}{L+A^2}dA(\frac{1}{\sqrt{L+A^2}}f_3,f_2)-\sqrt{L}(\frac{\sqrt{L}}{\sqrt{L+A^2}})^2 \frac{\sqrt{L}}{L+A^2}dA(f_2)-L(\frac{A}{\sqrt{L+A^2}})^2\\
&=&\frac{L}{(L+A^2)^2} d\alpha\wedge dA(f_3,f_2)-\frac{L^2}{(L+A^2)^2}dA(f_2)-\frac{L}{L+A^2}A^2.
\end{array}$$
Therefore
\begin{equation}\label{K}
K^{\infty}=\lim_{L\rightarrow \infty}K^L=-dA(f_2)-A^2.
\end{equation}

\section{The limit of Riemannian area element of $S$}\label{liRi}
It follows from (\ref{X1}) and (\ref{X2X3}) that 
$$X^1_L=\frac{\sqrt{L}}{\sqrt{L+A^2}}f^1,\mbox{ } X^2_L=f^2 \mbox{ and }X^3_L=\sqrt{L+A^2}f^3+\frac{A}{\sqrt{L+A^2}}f^1.$$ 
Therefore on $S$ we get
$$
d\sigma_L=X^2_L\wedge X^3_L=\sqrt{L+A^2}f^2\wedge f^3
$$
since that $f^1$ is null on $TS$. We can see that $\lim_{L\rightarrow \infty}K^Ld\sigma_L$ does not exist. In \cite{BTV}, to get an area form on $S$ it was necessary to multiply $d\sigma_L$ by $\frac{1}{\sqrt{L}}$ and take the limit as $L$ goes to infinity to obtain a surface form, i.e.,
$$
d\sigma=\lim_{L\rightarrow \infty }\frac{1}{\sqrt{L}}d\sigma_L=f^2\wedge f^3,
$$
which is  the Hausdorff measure on $S$. Therefore
\begin{equation}\label{KL}
\lim_{L\rightarrow \infty }\frac{1}{\sqrt{L}}K^Ld\sigma_L=K^{\infty}f^2\wedge f^3=K^{\infty}d\sigma.
\end{equation}

\section{The limit of normal curvatures of transverse curves in $S$}\label{curvnormal}

Suppose $\gamma(t)$ is a curve in $S$ such that $\gamma'(t)=a(t)f_2(\gamma(t))+b(t)f_3(\gamma(t))$, where $b(t)\neq 0$ for every $t$. Then $\gamma'(t)=a(t)X_2^L+b(t)\sqrt{L+A^2}X_3^L$ and the unitary tangent vector in the metric $g_L$ is
$$
T^L(t)=\frac{1}{\sqrt{a^2+b^2(L+A^2)}}(aX_2^L+b\sqrt{L+A^2}X_3^L).
$$
Let's write $a^L=\frac{a}{\sqrt{a^2+b^2(L+A^2)}}$ and $b^L=\frac{b\sqrt{L+A^2}}{\sqrt{a^2+b^2(L+A^2)}}$, so that 
$T^L=a^LX_2^L+b^LX_3^L$. Then
$$
\begin{array}{rcl}
\nabla^L_{T^L}T^L&=&\frac{d}{dt}a^LX_2^L+\frac{d}{dt}b^LX_3^L+a^L\nabla^L_{T^L}X_2^L+b^L\nabla^L_{T^L}X_3^L\vspace{.1cm}\\
&=&\frac{d}{dt}a^LX_2^L+\frac{d}{dt}b^LX_3^L+(-\sin\beta d\alpha(T^L)+\frac{\sqrt{L}}{2}\sin(2\beta)X^3_L(T^L))(a_LX_3^L-b^LX_2^L)\\
\end{array}
$$
Let be $N^L=-b^LX_2^L+a^LX_3^L$. Then $k_n^L=<\nabla^L_{T^L}T^L,N^L>$, so
$$
\begin{array}{rcl}
k_n^L&=&-b_L\frac{d}{dt}a^L+a_L\frac{d}{dt}b^L+(-\sin\beta d\alpha(T^L)+\frac{\sqrt{L}}{2}\sin(2\beta)X^3_L(T^L))\vspace{.1cm}\\
&=&\frac{1}{\sqrt{L+A^2}(a^2+b^2(L+A^2))}(abA\frac{d}{dt}A+(a\frac{d}{dt}b-b\frac{d}{dt}a)(L+A^2))\vspace{.1cm}\\
&&-\sin\beta(a_Ld\alpha(X_2^L)+b_Ld\alpha(X_3^L))+\sqrt{L}\sin\beta\cos\beta b_L\vspace{.1cm}\\
&=&\frac{1}{\sqrt{L+A^2}(a^2+b^2(L+A^2))}(abA\frac{d}{dt}A+(a\frac{d}{dt}b-b\frac{d}{dt}a)(L+A^2))\vspace{.1cm}\\
&&-\frac{A}{\sqrt{L+A^2}}\frac{a}{\sqrt{a^2+b^2(L+A^2)}}d\alpha(f_2)-\frac{A}{\sqrt{L+A^2}}\frac{b\sqrt{L+A^2}}{\sqrt{a^2+b^2(L+A^2)}}d\alpha(\frac{1}{\sqrt{L+A^2}}f_3)\vspace{.1cm}\\
&&+\sqrt{L}\frac{A}{\sqrt{L+A^2}}\frac{\sqrt{L}}{\sqrt{L+A^2}}\frac{b\sqrt{L+A^2}}{\sqrt{a^2+b^2(L+A^2)}}
\end{array}
$$
It follows from this formula that 
$$
k_n=\lim_{L\rightarrow\infty} k_n^L=A\frac{b}{|b|}
$$
We can see, from expressions $K^{\infty}=-dA(f_2)-A^2$ and $k_n=A\frac{b}{|b|}$, that both formulas depend only on $A$. This means that 
the "horizontal" geometry of the surface $S$ disappears  as $L$ goes to infinity. This will became clear as we proceed to find the rotation surfaces with $K^\infty$ constant, which are composed only of three families.

\section{The limit of length elements}

The length element in the metric $g^L$ on $\gamma$ is
$$
ds_L=
a^LX^2_L+b^LX^3_L.
$$
As $\lim _{L\rightarrow\infty}a^L=0$ and $\lim _{L\rightarrow\infty}b^L=\frac{b}{|b|}$, we get
$$
\lim _{L\rightarrow\infty}ds_L=\lim _{L\rightarrow\infty}(a^Lf^2+b^L\sqrt{L+A^2}f^3)=\frac{b}{|b|}f^3\lim _{L\rightarrow\infty}\sqrt{L+A^2}
$$
that does not exist. But as  in \cite{BTV} and section \ref{liRi}, if we multiply by $\frac{1}{\sqrt{L}}$ we obtain
$$
ds=\lim _{L\rightarrow\infty}\frac{1}{\sqrt{L}}ds_L=\frac{b}{|b|}f^3,
$$
which is Hausdorff measure for transversal curves. It follows that
\begin{equation}\label{kL}
\lim _{L\rightarrow\infty}\frac{1}{\sqrt{L}}k_n^Lds_L=k_nds=Af^3.
\end{equation}

\section{The Gauss-Bonnet theorem}
The proof of Gauss-Bonnet theorem in \cite{BTV} was done taking limits of  Gauss-Bonnet formulas in $(R^3,g_L)$ as $L$ goes to infinity:
$$
\int_S\frac{1}{ \sqrt{L}}K^L d\sigma_L+\int_{\partial S}\frac{1}{ \sqrt{L}}k_n^L ds_L=\frac{1}{ \sqrt{L}}2\pi\chi(S)
$$
to obtain $\int_SK^\infty d\sigma+\int_{\partial S}k_n ds=0$.

We will give below a straightforward proof due to  expressions of $K^{\infty}$ and $k_n$ obtained in sections \ref{curvgauss} and \ref{curvnormal}. We will restrict our theorem to regions where points are non singular and the boundary is constituted by transverse curves.

Let be $R\subset S$ a fundamental set, and $c$ a fundamental $2$-chain such that $|c|=R$. The oriented curve $\gamma=\partial c$ is the bounding curve of $R$. The curve $\gamma$ is piecewise differentiable, and composed by differentiable curves $\gamma_j:[s_j,s_{j+1}]\rightarrow S$, $j=1,\ldots,r$, with $\gamma_1(s_1)=\gamma_r(s_{r+1})$  and $\gamma_j(s_{j+1})=\gamma_{j+1}(s_{j+1})$, for $j=1,\ldots,r-1$. 
\begin{thm}(Gauss-Bonnet formula) Let ${R}$ be contained in a coordinate domain $U$ of $S$, let the bounding curve $\gamma$ of $R$ be a simple closed transverse curve. Then
$$\int_{\gamma}k_n+\int_R K^\infty=0,$$
where $k_n=\lim_{L\rightarrow\infty} k_n^L$  on $\gamma$ and $K^{\infty}=\lim_{L\rightarrow \infty}K^L$  on $R$.
\end{thm}
\Pf 
 From Stokes theorem and using (\ref{dAf3}) we get
$$
\int_RK^\infty=\int_cK^\infty f^2\wedge f^3=\int_c( -i(f_2)dA-A^2) f^2\wedge f^3=\int_c d(-Af^2)=-\int_{\partial c}Af^2=-\int_{\gamma}k_n.
$$


\section{The curvature $K^\infty$ for surfaces invariant by rotations}\label{rot}

Suppose $(a(t),b(t))$ is a curve in $\mathbb R^2$ and $(a(t),b(t),c(t))$ a horizontal curve in $\mathbb H^1$ such that 
$$
c(t)=\frac 1 2\int_0^t[a(s)b'(s)-b(s)a'(s)]ds.
$$
Rotating this curve we obtain a surface 
 $S$  \emph{invariant by rotations} given by 
$$f(u,v)=(a(v)\cos u-b(v)\sin u, b(v)\cos u+a(v)\sin u, c(v)).$$
The coordinate vector fields are 
$$
f_u=f_*(\frac{\partial}{\partial u})=(-a(v)\sin u-b(v)\cos u) e_1+(-b(v)\sin u+a(v)\cos u) e_2-\frac 1 2(a(v)^2+b(v)^2) e_3
$$
and
$$
f_v=f_*(\frac{\partial}{\partial v})=(a'(v)\cos u-b'(v)\sin u) e_1+ (a'(v)\sin u+b'(v)\cos u)  e_2.
$$
Then $f_v\in TS\cap D$, and as $<f_v,f_v>=a'(v)^2+b'(v)^2$ we obtain
\begin{equation}\label{f2}
f_2(u,v)=\frac{1}{\sqrt{a'(v)^2+b'(v)^2}}f_v(u,v).
\end{equation}
Therefore
$$
f_2(u,v)=\frac{a'(v)\cos u-b'(v)\sin u}{\sqrt{a'(v)^2+b'(v)^2}}e_1+ \frac{a'(v)\sin u+b'(v)\cos u}{\sqrt{a'(v)^2+b'(v)^2}} e_2 .
$$
Then
$$
\cos\alpha(u,v)= \frac{a'(v)\sin u+b'(v)\cos u}{\sqrt{a'(v)^2+b'(v)^2}},
$$
$$
\sin\alpha(u,v)=\frac{-a'(v)\cos u+b'(v)\sin u}{\sqrt{a'(v)^2+b'(v)^2}}
$$
and
$$
f_1(u,v)= \frac{a'(v)\sin u+b'(v)\cos u}{\sqrt{a'(v)^2+b'(v)^2}}e_1- \frac{a'(v)\cos u-b'(v)\sin u}{\sqrt{a'(v)^2+b'(v)^2}}e_2.
$$
As $<f_u(u,v),f_v(u,v)>=a(v)b'(v)-b(v)a'(v)$, we get 
$$
f_u-\frac{a(v)b'(v)-b(v)a'(v)}{\sqrt{a'(v)^2+b'(v)^2}}f_2=-\frac{a(v)a'(v)+b(v)b'(v)}{\sqrt{a'(v)^2+b'(v)^2}}f_1-\frac{a(v)^2+b(v)^2}{2}e_3
$$
so
$$
f_u-\frac{ab'-ba'}{\sqrt{(a')^2+(b')^2}}f_2=-\frac{a^2+b^2}{2}\left(e_3+2\frac{aa'+bb'}{\sqrt{(a')^2+(b')^2}(a^2+b^2)}f_1\right).
$$
It follows that
$$
f_u=\frac{ab'-ba'}{\sqrt{(a')^2+(b')^2}}f_2-\frac{a^2+b^2}{2}f_3
$$
where  $f_3=e_3+Af_1$ and 
$$
A=2\frac{aa'+bb'}{\sqrt{(a')^2+(b')^2}(a^2+b^2)}.
$$
From equation $d\sin\alpha=\cos\alpha\, d\alpha$ we get
$$
d\sin\alpha=\frac{a'\sin u+b'\cos u}{\sqrt{(a')^2+(b')^2}}(du+\frac{a''b'-b''a'}{(a')^2+(b')^2}dv)
$$
so
$$
d\alpha=du+\frac{a''b'-b''a'}{(a')^2+(b')^2}dv.
$$
Introducing polar coordinates 
$$a(t)=r(t)\cos\theta(t), \mbox{ }b(t)=r(t)\sin\theta(t)$$ 
we get $a^2+b^2=r^2$, $ab'-ba'=r^2\theta'$, $aa'+bb'=r r'$, $(a')^2+(b')^2=(r')^2+r^2(\theta')^2$, $a'b''-b'a''=(r^2(\theta')^2+2(r')^2-r r'')\theta'+r r' \theta''$. Taking a parameterization such that
$$
(a')^2+(b')^2=1,
$$
then
\begin{equation}\label{eqr}
A=2\frac{aa'+bb'}{(a^2+b^2)}=\frac{2 r r'}{r^2}=\frac{d}{dv}\ln r^2
\end{equation}
Also from (\ref{f2})
$$
f_2(u,v)=f_v(u,v).
$$
and (\ref{K})
\begin{equation}\label{kinf}
K^{\infty}=-dA(f_v)-A^2=-\frac{dA}{dv}-A^2
\end{equation}

\section{Surfaces of constant $K^\infty$  curvature invariant by rotations}\label{sol}
Suppose that $K^\infty$ is constant on $S$. Solving equation (\ref{kinf}) we obtain
$$
A(v)=-\sqrt{K^\infty}\tan(\sqrt{K^\infty} x+c_1),
$$
for $K^\infty>0$;
$$
A(v)=\sqrt{-K^\infty}\tanh(\sqrt{-K^\infty} x+c_1),
$$
for $K^\infty<0$;
$$
A(v)=\frac{1}{x+c_1}
$$
for $K^\infty=0$.
From (\ref{eqr}) we get:
$$
r(v)=c_2\sqrt{\cos(\sqrt{K^\infty}v+c_1)}
$$
for $K^\infty>0$;
$$
r(v)=c_2\sqrt{\cosh(\sqrt{-K^\infty}v+c_1)}
$$
for $K^\infty<0$;
$$
r(v)=c_2\sqrt{v+c_1}
$$
for $K^\infty=0$.
In any case 
\begin{equation}\label{tc}
\begin{array}{rcl}
\theta(v)&=&\int_0^v\frac{\sqrt{1-r'(t)^2}}{r(t)}dt\vspace{0.1cm}\\
c(v)&=&\frac{1}{2}\int_0^vr(t)\sqrt{1-r'(t)^2}dt.
\end{array}
\end{equation}

\subsection{Case $K^\infty>0$}
We can write
$$
r(v)=r_0\sqrt{\cos(\sqrt{K^\infty}v)},
$$
and $c(v)$, $\theta(v)$ as in (\ref{tc}). The surface is defined only for 
$|v|<\frac{1}{\sqrt{K^\infty}}\cos^{-1}(-\frac{2}{r_0^2 K^\infty}+\sqrt{\frac{4}{r_0^4 (K^\infty)^2}+1})$
 and,  in figure \ref{positivo} below, we can see the graphics of $S$ with $r_0=1$ and $K^\infty=1$, indicating the horizontal curve.

\begin{figure}[h]
\center
\includegraphics[width=.5\textwidth, angle=0]{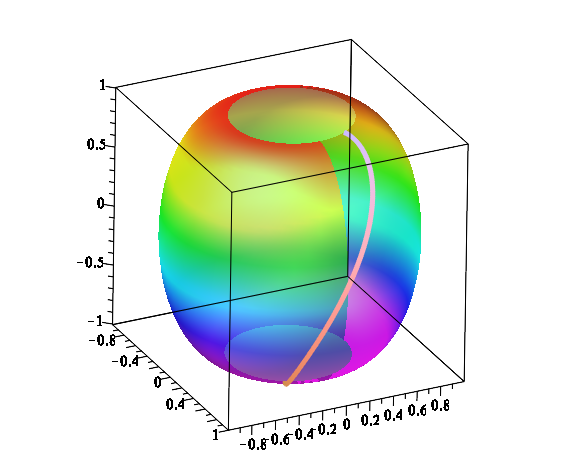}
\caption{Surface with $r_0=1$ and $K=1$  }
\label{positivo}
\end{figure}

\subsection{Case $K^\infty=0$}
We can write
$$
r(v)=r_0\sqrt{v},
$$
and $c(v)$, $\theta(v)$ as in (\ref{tc}). The curve is defined only for $v>\frac{R^2}{4}$
and,  in figure \ref{nulo} below, we can see the graphics of $S$ with $r_0=1$ and $K^\infty=0$, showing the horizontal curve.
\begin{figure}[h]
\center
\includegraphics[width=.4\textwidth, angle=0]{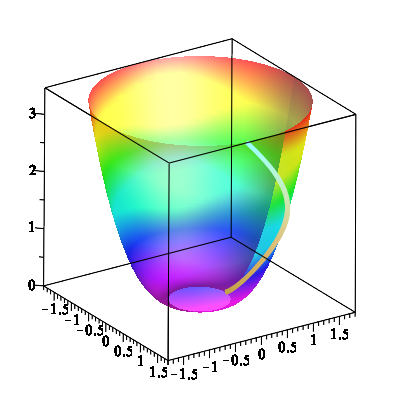}
\caption{Surface with $r_0=1$ and $K^\infty=0$  }
\label{nulo}
\end{figure}

\subsection{Case $K^\infty<0$}
We can write
$$
r(v)=r_0\sqrt{\cosh(\sqrt{-K^\infty}v)}
,
$$
and $c(v)$, $\theta(v)$ as in (\ref{tc}). The surface is defined only for 
$|v|<\frac{1}{\sqrt{-K^\infty}}\cosh^{-1}(\frac{-2}{r_0^2K^\infty}+\sqrt{1+\frac{4}{r_0^4(K^\infty)^2}})$ 
and, in figure \ref{negativo} below, we can see the graphics of $S$ with $r_0=1$ and $K^\infty=-1$, showing the horizontal curve.
\begin{figure}[h]
\center
\includegraphics[width=.4\textwidth, angle=0]{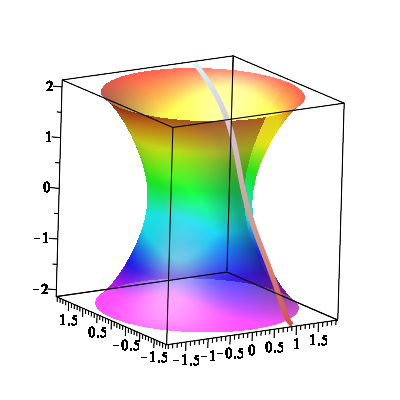}
\caption{Surface with $r_0=1$ and $K^\infty=-1$  }
\label{negativo}
\end{figure}

\end{document}